\documentclass{article}
\usepackage{amsmath}
\usepackage{amsfonts}
\usepackage{hyperref}
\usepackage{graphicx}
\usepackage[indention=10pt]{subcaption} 
\usepackage{tikz}
\usetikzlibrary{arrows, shapes, calc}

\tikzset{
  base/.style = {circle, align=center, inner sep=0pt, text centered, font=\sffamily, black, draw=black, text width=1.5em, very thick},
  inf/.style = {red, draw=red},
  opt/.style = {black!30!green, draw=black!30!green},
  bound/.style = {black!50!red, draw=black!50!red},
  grayed/.style = {black!30, draw=black!30},
  xlabel/.style = {rectangle, text width=16mm},
  zlabel/.style = {rectangle, minimum height=7mm, text width=16mm},
  edge/.style = {above left, draw=none},
  offshoot/.style = {cloud, draw, cloud puffs=9.178, cloud puff arc=120, aspect=2, fill=white, text width=1.5cm}
}

\numberwithin{equation}{section}
\newtheorem{example}{Example}

\newcommand{\set}[1]{\left\{#1\right\}}

\begin{document}

\title{Shaping and Trimming Branch-and-bound Trees\footnote{This paper is a written version of the talk given by the first author at the MIP 2017 Workshop in Montreal, Canada.}}
\author{Philipp M. Christophel\thanks{philipp.christophel@sas.com, SAS Institute, Inc.}\and Imre P\'olik\thanks{imre@polik.net, SAS Institute, Inc.}}

\date{\today}
\maketitle

\begin{abstract}
We present a new branch-and-bound type search method for mixed integer linear optimization problems based on the concept of \emph{offshoots} (introduced in this paper). While similar to a classic branch-and-bound method, it allows for changing the order of the variables in a dive (\emph{shaping}) and removing unnecessary branching variables from a dive (\emph{trimming}).  The regular branch-and-bound algorithm can be seen as a special case of our new method. We also discuss extensions to our new method such as choosing to branch from the top or the bottom of an offshoot. We present several numerical experiments to give a first impression of the potential of our new method.
\end{abstract}

\section{Introduction}


In this paper we are discussing mixed integer linear optimization problems (MILP), i.e., optimization problems with a linear objective function, linear constraints, and integrality restrictions on some or all of the variables. Most of the techniques described in the paper generalize naturally to problems with non-linear constraints as well (MINLP). The typical approach to solve any optimization problem with an integrality restriction on the feasible domain involves variants of the branch-and-bound algorithm first described for general integer optimization by Land and Doig in \cite{Land1960}. For details about the origins of branch-and-bound, see also \cite{Cook2012}. The branch-and-bound method is at the core of every software to solve mixed integer optimization problems and is successfully used to solve a variety of practical problems. But it is also known that in the worst case the branch-and-bound method will enumerate all possible solutions, leading to a disastrous performance. In Example \ref{badtree} we show such a case.

\begin{example}\label{badtree}
Consider the following integer optimization problem with three binary variables:
\begin{align}
\notag \min\  x_1 - 2 x_2 - 6 x_3 \\
\notag         -3x_1 - 4 x_2 - 2 x_3 &\ge -8 \\
\label{eq:badtree}        3x_1 - 4 x_2 - 2 x_3 &\ge -5 \\
\notag         -3x_1 + 4 x_2 - 2 x_3 &\ge -4 \\
\notag         3x_1 + 4 x_2 - 2 x_3 &\ge -1 \\
\notag        x \in \{0, 1\}^3
\end{align}
We solve this example with a depth-first branch-and-bound algorithm where in the left nodes variables are fixed to one. In this example, it does not matter which branching variable selection method is used since there is only one fractional variable in each node and a traditional branch-and-bound method chooses the branching variable only from the fractional variables. In Figure \ref{fig:badtree} we show the branch-and-bound tree resulting from this example. We show the objective value of the LP relaxation below each node and the solution values for the three variables to the right of each node. The number inside the nodes shows the order in which they are processed.
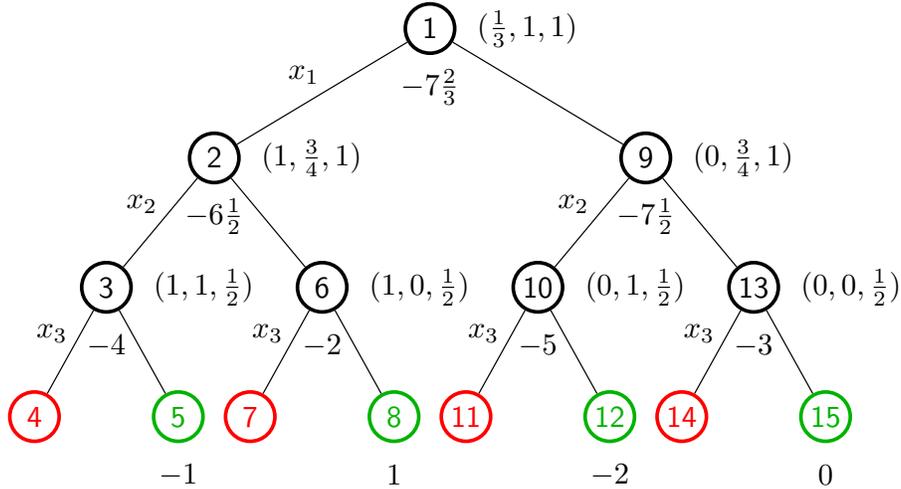
\begin{figure}[h!]
\resizebox{\textwidth}{!}{%
\begin{tikzpicture}[-,>=stealth',
level/.style={sibling distance = 5cm/#1, level distance = 1.5cm},
every node/.style={base}
] 
\node [label={[zlabel]below:{$-7\frac{2}{3}$}}, label={[xlabel]right:{$(\frac{1}{3}, 1, 1)$}}] {1}
    child{ node [label={[zlabel]below:{$-6\frac{1}{2}$}}, label={[xlabel]right:{$(1, \frac{3}{4}, 1)$}}] {2}
            child{ node [label={[zlabel]below:{$-4$}}, label={[xlabel]right:{$(1, 1, \frac{1}{2})$}}] {3}
            	child{ node [inf] {4} edge from parent node[edge] {$x_3$}}
		child{ node [opt, label={[zlabel]below:{$-1$}}] {5}}
		edge from parent node[edge] {$x_2$}
            }
            child{ node [label={[zlabel]below:{$-2$}},label={[xlabel]right:{$(1, 0, \frac{1}{2})$}}] {6}
              	child{ node [inf] {7} edge from parent node[edge] {$x_3$}}
		child{ node [opt, label={[zlabel]below:{$1$}}]{8}}
            }
            edge from parent node[edge] {$x_1$}                    
    }
    child{ node [label={[zlabel]below:{$-7\frac{1}{2}$}}, label={[xlabel]right:{$(0, \frac{3}{4}, 1)$}}] {9}
            child{ node [label={[zlabel]below:{$-5$}}, label={[xlabel]right:{$(0, 1, \frac{1}{2})$}}] {10} 
		child{ node [inf] {11} edge from parent node[edge] {$x_3$}}
		child{ node [opt, label={[zlabel]below:{$-2$}}] {12}}
		edge from parent node[edge] {$x_2$}
            }
            child{ node [label={[zlabel]below:{$-3$}}, label={[xlabel]right:{$(0, 0, \frac{1}{2})$}}] {13}
		child{ node [inf] {14} edge from parent node[edge] {$x_3$}}
		child{ node [opt, label={[zlabel]below:{$0$}}] {15}}
            }
    }
;
\end{tikzpicture}}
\caption{The branch-and-bound tree for the optimization problem \ref{eq:badtree}.}\label{fig:badtree}
\end{figure}

Note that the optimal objective value is $-2$ but it is obtained only after 12 nodes have been processed and proving its optimality requires $15$ nodes, the maximal number of nodes possible, which corresponds to enumerating all possible solutions, i.e., the leaf nodes of the tree in Figure \ref{fig:badtree}. Note that the leaf nodes are alternating between infeasible and integral nodes. This is the reason why we have such a bad tree in this example; the enumeration ended up with different types of leaf nodes next to each other. Hence the branch-and-bound method has no chance to prune several leaf nodes together early on at a common ancestor.
\end{example}

Research on branch-and-bound algorithms has put a huge emphasis on making the decisions in the algorithm in such a way as to avoid enumeration of large parts of the solution space. Especially selecting the branching variables has been studied extensively (see, for example, \cite{Achterberg2005}), because which branching variable is chosen determines the shape of the tree of a branch-and-bound method and thus how many nodes need to be processed. In this paper we will present an implementation for a branch-and-bound method that follows a different approach in which the shape of the tree can be changed. This means that we want to have a branch-and-bound algorithm where the decisions on the branching variables can be deferred to a time when potentially more information is available to make these decisions. Example \ref{goodtree} demonstrates that the same problem from Example \ref{badtree} can be solved more efficiently if the branch-and-bound tree has a different shape.

\begin{example}\label{goodtree}
In Figure \ref{fig:goodtree} we show the branch-and-bound tree for the same mixed integer optimization problem as in Example \ref{badtree} with the only difference that we reversed the order in which we branched on the variables, i.e., we branched on $x_3$ first although it is not even fractional. In this case we only need to process $9$ nodes; we show the remaining nodes (without a number) just to highlight the structure of the tree.

\begin{figure}[h!]
\resizebox{\textwidth}{!}{%
\centering
\begin{tikzpicture}[-,>=stealth',
level/.style={sibling distance = 5cm/#1, level distance = 1.5cm},
every node/.style={base}
] 
\node [label={[zlabel]below:{$-7\frac{2}{3}$}}, label={[xlabel]right:{$(\frac{1}{3}, 1, 1)$}}] {1}
    child{ node [label={[zlabel]below:{$-7\frac{2}{3}$}}, label={[xlabel]right:{$(\frac{1}{3}, 1, 1)$}}] {2}
            child{ node [label={[zlabel]below:{$-7\frac{2}{3}$}}, label={[xlabel]right:{$(\frac{1}{3}, 1, 1)$}}] {3}
            	child{ node [inf] {4} edge from parent node[edge] {$x_1$}}
		child{ node [inf] {5}}
		edge from parent node[edge] {$x_2$}     
            }
            child{ node [label={[zlabel]below:{$-5\frac{2}{3}$}}, label={[xlabel]right:{$(\frac{1}{3}, 0, 1)$}}] {6}
              	child{ node [inf] {7} edge from parent node[edge] {$x_1$}}
		child{ node [inf]{8}}
            }
            edge from parent node[edge] {$x_3$}                    
    }
    child{ node [opt, label={[zlabel]below:{$-2$}}, label={[xlabel]right:{$(0, 1, 0)$}}] {9}
            child{ node [grayed] {} 
		child{ node [opt] {} edge from parent node[edge] {$x_1$}}
		child{ node [opt] {}}
		edge from parent node[edge] {$x_2$}
            }
            child{ node [grayed] {}
		child{ node [opt] {} edge from parent node[edge] {$x_1$}}
		child{ node [opt] {}}
            }
    }
;
\end{tikzpicture}}
\caption{Problem \eqref{eq:badtree} solved with a different branch-and-bound tree.}\label{fig:goodtree}
\end{figure}
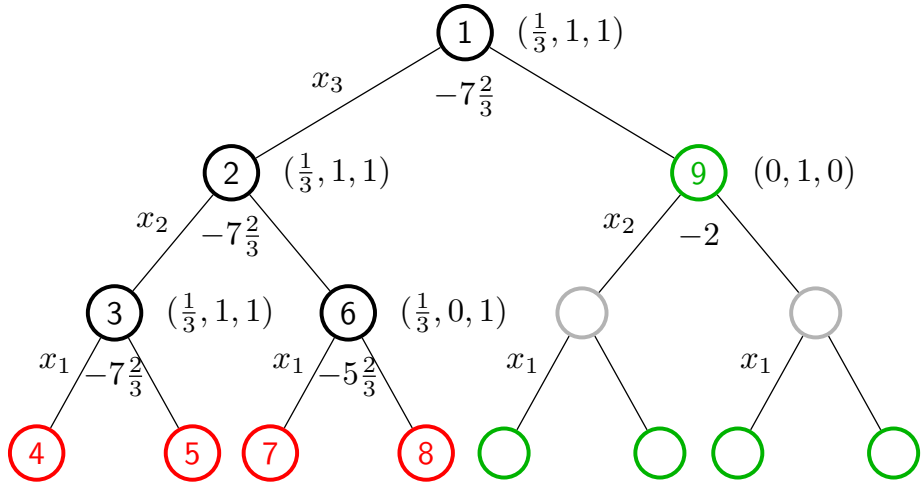

Note that now the infeasible leaf nodes are all on the left side of the tree and all the integral nodes are on the right side of the tree. The integral leaf nodes do not all have to be processed because the common ancestor, node $9$, covers them all. Also, if we can detect the integer infeasibility of node 2 (for example with probing or some other node presolver technique), then we can prune all the nodes below it and solve the problem even faster. The lesson is that it is better to have a tree in which nodes with similar properties are the leaf nodes of subtrees so that they can be dealt with at a common ancestor.
\end{example}

The question here is twofold. First, we need to reshape the tree if we believe that the current tree structure is inefficient. The second question is equally important: we need to do this in an efficient way, so that we can preserve most of the advanced techniques that make branch-and-bound implementations perform well in practice. While it is possible to simply throw away most of the tree and roll back to the last known good decision point (or variants thereof called restarting as, for example, discussed in \cite{Achterberg2005}), we want to explicitly look into other possibilities here.

There has been some research in branch-and-bound methods with an adjustable (sometimes called dynamic) tree. The earliest we are aware of is by Glover and Tangedahl \cite{Glover1976}. Chv{\'a}tal in \cite{chvatal1997} and then Hanafi and Glover in \cite{Hanafi2002} revisited the topic. These papers give valuable insights into alternative methods for solving mixed integer optimization problems but unfortunately do not discuss the implementational challenges. Furthermore, resolution search from \cite{chvatal1997}, for example, is not similar enough to a classic branch-and-bound method such that many of the methods modern solvers successfully use to solve problems today are not directly applicable.

\section{Diving, Shaping, and Trimming}
In this section we discuss three very important concepts for the remainder of this paper: diving, shaping, and trimming. We do so using a depth-first branch-and-bound method because these concepts are easier to explain in this method and are also a natural extension to it.

Depth first branch-and-bound (sometimes also called last-in-first-out, i.e. LIFO, branch-and-bound) is a variant of branch-and-bound where the next node processed is always the most recent node added to a stack of open nodes. In practice it is possible to store the open nodes with a stack of bound changes. The depth-first branch-and-bound method also minimizes the number of open nodes. The result is a very memory-efficient branch-and-bound method. 

In the depth-first branch-and-bound method we repeatedly go down the tree only changing one variable at a time. We call this process of going down a tree \emph{diving}. Another advantage of the depth-first branch-and-bound method is that the LP relaxations during diving can be solved very efficiently using a dual simplex algorithm where most data structures (most importantly the factorization of the basis matrix) can be kept up to date. We call this \emph{hotstarting} the dual simplex to express that it is even better than \emph{warmstarting}, which typically implies that a known dual feasible basis is used to initialize the dual simplex algorithm. When backtracking in the depth-first branch-and-bound we cannot use hotstarting, but since the difference between nodes is typically small we can warmstart from the last basis instead of resolving from scratch. Since diving is much more efficient, current implementations of non-depth-first branch-and-bound methods also use it to process nodes quickly and only do a full node selection if the current dive does not seem promising anymore.

The disadvantage of the depth-first branch-and-bound method is that the problem described in the introduction is aggravated: a bad decision early on can result in a very bad enumeration tree and thus long running time or a failure to solve the problem within some resource limitation. But it is also much easier to revise earlier decisions and change the order of bound changes in a dive. Notice that for the status of the final node in a dive the order in which the variables were fixed does not matter. The order of the bound changes in a dive in some sense defines the shape of the branch-and-bound tree. Hence we call changing the order of the bound changes in a dive \emph{shaping}. Since in a depth-first branch-and-bound we store the bound changes in a stack we can decide to undo them in a different order than we did them during the dive. The only thing we have to keep in mind is that we can only change the order of the bound changes up to the last node where we have already explored the other side of the bound change.

In a depth-first branch-and-bound algorithm a dive has to end in a pruned node. A node is pruned either because the LP relaxation is infeasible or because the objective value exceeds the cutoff \footnote{The case of an integral solution can be seen as first establishing a new cutoff and then pruning the node.}. It is possible that a situation occurs where a dive contains more bound changes than are strictly necessary to prune a node. In this case it is possible to remove the unneeded bound changes from the dive before backtracking. Since this trims the dive down to a smaller set of bound changes we call this \emph{trimming}. 

There are a number of ways to trim a dive. For problems with general integer variables it is possible to remove multiple bound changes on the same bound of the same variable and keep only the tightest one. It is also possible to use reduced cost or Farkas certificate values to trim dives. In fact, this problem is identical to the one we are facing when trying to identify an irreducible infeasible system (IIS), so all the reduction techniques in that domain apply readily to our setup; see \cite{chinneck2008} for details. In the following sections we will sample a few methods.

Shaping and trimming clearly can improve a depth-first branch-and-bound implementation a lot, and the implementational complexity is very low. For shaping, the obvious difficulty is to come up with good rules on which bound change should be undone first. But our experience has been that even simple rules already lead to an improvement. For trimming, the trade-off is between time spent trimming the tree and simply processing nodes. But here as well simple strategies already yielded benefits so that it should be possible to improve any depth-first branch-and-bound implementation not making use of trimming significantly.

The only downside is that if node presolving techniques are used in a depth-first branch-and-bound method it is necessary to keep track of implied bound changes separately from the actual branching decisions. As a result, during backtracking some tightenings from node presolve have to be redone.

The concepts of shaping and trimming the tree already appear in principle in \cite{Glover1976}, but that paper does not include any implementational considerations.

\section{A New Branch-and-bound Method}
\begin{figure}[h!]
\centering
\begin{tikzpicture}[-,>=stealth',
every node/.style={base},
level/.style={ level distance = 2cm},
] 
\node [ label={right:$F$}, label={below:$z^*$}] {s}
    child{ node [cloud, draw, cloud puffs=9.178, cloud puff arc=120, aspect=2, text width=1.5cm] {$D$}
            	child{ node [inf]{t} 
            }
            child{edge from parent[draw=none]}
    }
    child {edge from parent[draw=none]}
;
\end{tikzpicture}
\caption{The structure of an offshoot.\label{offshoot}}
\end{figure}
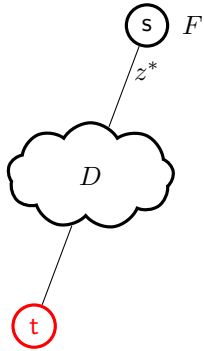

In this section we present the  basic idea of a new branch-and-bound method that allows for shaping and trimming but is not a depth-first method. The fundamental idea is to perform a branch-and-bound method on objects we call \emph{offshoots} instead of performing it on individual nodes. An offshoot (see  Figure~\ref{offshoot}) is an object that represents a collection of nodes in a tree. It consists of a top node $s$, represented by a set $F$ of initial bound changes, with an attached set $D$ of bound changes representing a dive in the branch-and-bound tree. Applying both the initial bound changes in $F$ and the bound changes in the dive $D$ has to result in a node $t$ that can be pruned, either because it is infeasible or because its objective value exceeds the current cutoff.\footnote{The cutoff is derived from the currently best known primal feasible solution.} Note that the order of the bound changes in $D$ is not determined,\footnote{In a practical implementation we can remember the original order of bound changes so that we can use the intermediate objective values to prune undisturbed nodes inside an offshoot.} only the set of all the bound changes needed to reach a terminal node.

Instead of storing a set of open nodes that still need to be processed we store a set of open offshoots. An offshoot is considered open if it has bound changes in its dive that have not been processed. Once the list of open offshoots is empty the problem is solved.

This new method begins with creating a first offshoot for which the set $F_0$ of initial bound changes is empty. Then it performs a dive until it reaches a node that can be pruned and stores the bound changes of this first dive in the set $D_0$ of the first offshoot. Then the first offshoot is added to the list of open offshoots.

From now on, in each iteration, the method selects an offshoot from the list of open offshoots as parent offshoot $p$ for a new offshoot $k$ to create. The method also needs to select a bound change to process associated with an offshoot variable $i$ from the list $D_p$ of unprocessed dive bound changes of its parent. The initial set of bound changes for the new offshoot $k$ is $F_k = F_p \cup (D_k \setminus \set{x_i \le b}) \cup \set{x_i \ge b+1}$ if the bound change for the selected variable was branching down or $F_k = F_p \cup ( D_k \setminus \set{x_i \ge b}) \cup \set{x_i \le b-1}$ if it was branching up. The new offshoot starts with a node that corresponds to a right node of the dive but since we can freely choose from all bound changes in the dive it might be a right node that does not correspond to any of the dive nodes that were processed when the offshoot was created. This choosing of the variable from the dive corresponds to shaping the tree.

After creating the initial node of the new offshoot we solve the LP relaxation of the top node in the new offshoot. If the top node can be pruned we proceed by selecting a new parent offshoot right away. Otherwise we store the objective value as the top bound $z_k^*$ of the new offshoot. Then we perform a dive until we reach a node that can be pruned either because it is infeasible or because its objective value exceeds the current cutoff. If we encounter a new primal feasible solution we update the cutoff. When updating the cutoff we can also remove all open offshoots for which the top bound exceeds the cutoff.\footnote{In addition, we can also remove those nodes inside offshoots that have not been disturbed yet if their objective value exceeds the cutoff.}

In this setup we can also easily perform trimming. As mentioned before, this can be done, for example, by removing multiple bound changes on the same bound of a variable (only in the case of general integer variables) or by inspecting the dual information of the pruned node. To specify in more detail: the dual information vector $r$ is either the reduced cost vector for cutoff nodes or the Farkas certificate for infeasible nodes. An upper bound change on variable $i$ can be removed if $r_i \ge 0$, and a lower bound change on variable $i$ can be removed if $r_i \le 0$. 

After trimming the dive we can store the new offshoot in the list of open offshoots and remove the parent offshoot if all the bound changes in its dive have been processed.

\begin{figure}[t!]
\begin{subfigure}[t]{0.5\textwidth}
\centering
\resizebox{\textwidth}{!}{%
\begin{tikzpicture}[-,>=stealth',
level/.style={sibling distance = 5cm/#1, level distance = 1.5cm},
every node/.style={base}
] 
\node [label={[xlabel]right:{$(\frac{1}{3}, 1, 1)$}}] {1}
    child{ node [label={[xlabel]right:{$(1, \frac{3}{4}, 1)$}}] {2}
            child{ node  [label={[xlabel]right:{$(1, 1, \frac{1}{2})$}}] {3}
            	child{ node [inf] {4} edge from parent node[edge] {$x_3$}}
		child{edge from parent[draw=none]}
		edge from parent node[edge] {$x_2$}     
            }
            child{edge from parent[draw=none]}
            edge from parent node[edge] {$x_1$} 
    }
    child{edge from parent[draw=none]}
;
\end{tikzpicture}}
\caption{The initial dive \dots}
\end{subfigure}
\begin{subfigure}[t]{0.5\textwidth}
\centering
\resizebox{\textwidth}{!}{%
\begin{tikzpicture}[-,>=stealth',
level/.style={sibling distance = 5cm/#1, level distance = 1.5cm},
every node/.style={base}
] 
\node {1}
    child{
            child{ 
            	child{ node [inf] {4}}
		child{edge from parent[draw=none]}
		edge from parent node[offshoot] {$x_1 \ge 1$ \\ $x_2 \ge 1$ \\ $x_3 \ge 1$}
            }
            child{edge from parent[draw=none]}
    }
    child{edge from parent[draw=none]}
;
\end{tikzpicture}}
\caption{\dots becomes the first offshoot.}
\end{subfigure}
\begin{subfigure}[t]{0.5\textwidth}
\centering
\resizebox{\textwidth}{!}{%
\begin{tikzpicture}[-,>=stealth',
level/.style={sibling distance = 5cm/#1, level distance = 1.5cm},
every node/.style={base}
] 
\node {1}
    child{  node[offshoot] {$x_2 \ge 1$ \\ $x_3 \ge 1$}
           child {node [grayed] {}
	          		child{ node [inf] {4} edge from parent node[edge] {$x_1$}}
			child{ node [inf] {5}}
	}
	child{edge from parent[draw=none]}
    }
    child{edge from parent[draw=none]}
;
\end{tikzpicture}}
\caption{We choose $x_1$ as the first offshoot variable but the top node of the new offshoot is immediately infeasible.}
\end{subfigure}
\begin{subfigure}[t]{0.5\textwidth}
\centering
\resizebox{\textwidth}{!}{%
\begin{tikzpicture}[-,>=stealth',
level/.style={sibling distance = 5cm/#1, level distance = 1.5cm},
every node/.style={base}
] 
\node [label={[zlabel]below:{$-7\frac{2}{3}$}}] {1}
    	child { node [grayed] {}
    		child {node [grayed] {}
	          		child{ node [inf] {4} edge from parent node[edge] {$x_1$}}
			child{ node [inf] {5}}
			edge from parent node[edge] {$x_2$}
		        }
  	          child{edge from parent[draw=none]}
	          child{node [label={[zlabel]below:{$-5\frac{2}{3}$}}, label={[xlabel]right:{$(\frac{1}{3}, 0, 1)$}}] {6}
	          		child {node[offshoot] {$x_1 \ge 1 $} 
	          		child{ node [inf] {7}}
	                     child{edge from parent[draw=none]}
	                     }
	                     child{edge from parent[draw=none]}
	          }
                     edge from parent node [offshoot, above left] {$x_3 \ge 1$}
	}
         child{edge from parent[draw=none]}
;
\end{tikzpicture}}
\caption{As the second offshoot variable we choose $x_2$ and create a second open offshoot with the top node labeled $6$.}
\end{subfigure}
\begin{subfigure}[t]{0.5\textwidth}
\centering
\resizebox{\textwidth}{!}{%
\begin{tikzpicture}[-,>=stealth',
level/.style={sibling distance = 5cm/#1, level distance = 1.5cm},
every node/.style={base}
] 
\node {1}
    	child { node [grayed] {}
    		child {node [grayed] {}
	          		child{ node [inf] {4} edge from parent node[edge] {$x_1$}}
			child{ node [inf] {5}}
			edge from parent node[edge] {$x_2$}
		        }
  	          child{edge from parent[draw=none]}		        
	          child{node {6}
	          		child {node[offshoot] {$x_1 \ge 1 $} 
	          		child{ node [inf] {7}}
	                     child{edge from parent[draw=none]}
	                     }
	                     child{edge from parent[draw=none]}
	          }
	          edge from parent node[edge] {$x_3$}
	}
         child{node [opt, label={[zlabel]below:{$-2$}}, label={[xlabel]right:{$(0, 1, 0)$}}] {8}}
;
\end{tikzpicture}}
\caption{We have two offshoots to choose from. Since its top bound is better we choose the first offshoot and choose the last remaining bound change. The resulting top node of the new offshoot (8) is integer so we do not create a new offshoot.}
\end{subfigure}
\begin{subfigure}[t]{0.5\textwidth}
\centering
\resizebox{\textwidth}{!}{%
\begin{tikzpicture}[-,>=stealth',
level/.style={sibling distance = 5cm/#1, level distance = 1.5cm},
every node/.style={base}
] 
\node {1}
    	child { node [grayed] {}
    		child {node [grayed] {}
	          		child{ node [inf] {4} edge from parent node[edge] {$x_1$}}
			child{ node [inf] {5}}
			edge from parent node[edge] {$x_2$}
		        }
	          child{node {6}
	          	       child{ node [inf] {7}}
	                  child{ node [inf] {9}}
	          }
	          edge from parent node[edge] {$x_3$}
	}
         child{node [opt] {8}}
;
\end{tikzpicture}}
\caption{We choose the only bound change in the only open offshoot and the resulting offshoot is infeasible in the top node. The method stops with optimal objective value $-2$.}
\end{subfigure}
\caption{A step-by-step example of the new method. \label{example}}
\end{figure}
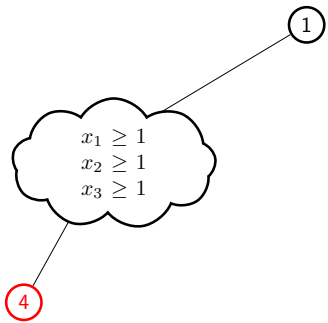
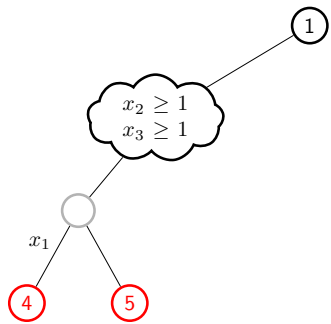
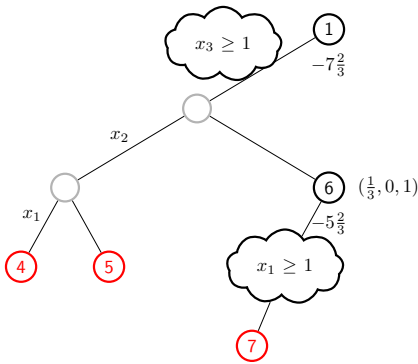
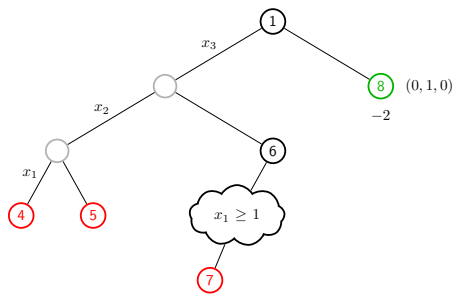
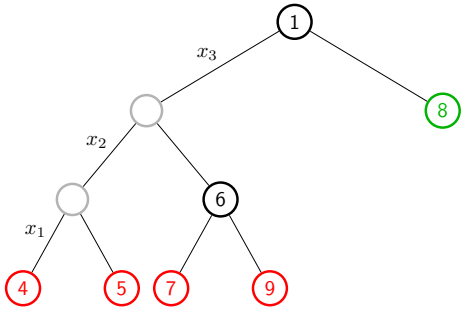

This continues until the list of open offshoots is empty. Figure \ref{example} shows an example where the new method is applied to Example \ref{badtree}.

\section{Improvements and Extensions}\label{improv}
As with many similar methods it is necessary to improve and extend our new method to get the best possible performance. In this section we list some more or less obvious ways to overcome some of the weaknesses of the new method.

\subsection{Branching From the Top}
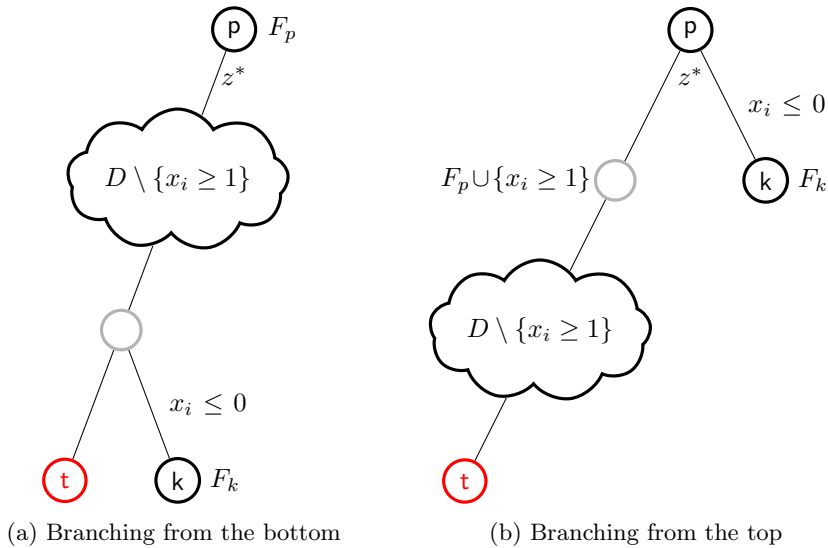
\begin{figure}
\begin{subfigure}[t]{0.5\textwidth}
\centering
\begin{tikzpicture}[-,>=stealth',
every node/.style={base},
level/.style={ level distance = 2cm}
] 
\node [ label={right:$F_p$}, label={below:$z^*$}] {p}
    child{ node [cloud, draw, cloud puffs=9.178, cloud puff arc=120, aspect=2, text width=2cm] {$D \setminus \set{x_i \ge 1}$}
            child{ node [grayed]{}
            child{ node [inf]{t} }
            child{ node [label={right:$F_k$}] {k}
            	edge from parent node[edge, right, text width=1.5cm] {$x_i \le 0$} 
            }
            }
            child {edge from parent[draw=none]}
    }
    child {edge from parent[draw=none]}
;
\end{tikzpicture}
\caption{Branching from the bottom}
\end{subfigure}
\begin{subfigure}[t]{0.5\textwidth}
\centering
\begin{tikzpicture}[-,>=stealth',
every node/.style={base},
level/.style={ level distance = 2cm, sibling distance = 2cm}
] 
\node [label={below:$z^*$}] {p}
    child { node [grayed,  label={left, text width=2cm:$F_p \cup \set{x_i \ge 1}$}] {}
    	child{ node [cloud, draw, cloud puffs=9.178, cloud puff arc=120, aspect=2, text width=2cm] {$D \setminus \set{x_i \ge 1}$}
            	child{ node [inf]{t}}
           	child {edge from parent[draw=none]}
	}
          child {edge from parent[draw=none]}
    }
    child{node [label={right:$F_k$}] {k}
            	edge from parent node[edge, right, text width=1.5cm] {$x_i \le 0$} 
    }
;
\end{tikzpicture}
\caption{Branching from the top}
\end{subfigure}
\caption{The two ways to branch illustrated.}\label{bottomtop}
\end{figure}

In the description of the method in the previous section we only added new offshoots below their parent. This can be seen as \emph{branching from the bottom} of an offshoot. It is also possible to \emph{branch from the top} of an offshoot. Then the new offshoot inherits only the bound changes its parent had at the top and additionally exactly one bound change from the dive flipped to the other side. The parent is then adjusted as well and one of the bound changes is moved from the dive to the top. To be precise, the initial set of bound changes for the new offshoot $k$ is $F_k = F_p \cup \set{x_i \ge b+1}$ if the bound change for the selected variable was branching down or $F_k = F_p \cup \set{x_i \le b-1}$ if it was branching up. Figure \ref{bottomtop} illustrates both types of branching next to each other.

The big advantage of this additional level of flexibility is that  we can decide which type of branching to use based on how sure we are that an offshoot variable is a good choice. If we are not sure whether a bound change will have large impact and hence should be at the top of the tree we can choose to branch from the bottom to minimize the effect if we made a bad choice. If, on the other hand, we have a strong indication that a bound change will have a huge impact and should be at the top of the tree, then we can branch from the top.

\subsection{Creating Offshoots and Advanced Trimming}\label{advtrim}
For the correctness of the method it is not necessary to create offshoots by diving. Any method that creates a set of bound changes that results in a pruned node can be used for the diving set in an offshoot.

One slight modification to the method is to apply several bound changes at once before solving an LP relaxation. We call this \emph{plunging}. This can go as far as fixing all integer variables since the result is guaranteed to be pruned and trimming can then be used to reduce the set of bound changes.

Another possibility is to use conflict analysis as described in \cite{Achterberg2005} to obtain a clause. For offshoots that end in an infeasible node the dive set of bound changes is precisely a clause. Hence it is also possible to apply the method described by Karzan  et al.~in \cite{Karzan2009} to obtain a minimal clause using a MIPing approach.

\subsection{Improved Pruning}
One of the disadvantages of the new method is that pruning by bound after a new primal feasible solution has been found is complicated. Obviously, we can prune whole offshoots as soon as their top bound exceeds the new bound. But it can happen that for some offshoots the top bound is not large enough although applying some of the bound changes from the dive would result in an LP bound that would lead to a pruning.

This issue can be overcome partially by storing the objective values obtained during a dive. As long as no new offshoot is created from the dive (or new offshoots are created only from the bottom and in the order of the original dive) we can use the objective values to trim the dives after a new bound has been found. Since trimming also invalidates the bounds from the dive it is advisable to delay trimming until we first want to create an offshoot. This requires slightly more memory since more bound changes and dual information might have to be stored, but it could result in significantly better performance.

\subsection{Bounding Offshoots}

When branching from the top, the top bound of an offshoot remains a valid bound on all the nodes below this offshoot. But since we add a bound change to the top of the offshoot the bound is obviously not as strong as it could be. Hence it might be worthwhile updating the top bound after branching from the top.

We propose three methods of increasing computational effort to strengthen the bound. The first method is to derive a bound on the top node of the parent offshoot by using the reduced cost of the just-solved top node of the new offshoot. 

The second method is a bit more general but also requires more computational effort. It involves simply evaluating the dual solution of the new offshoot's top node for the bounds of the parent offshoot.

The third method is to solve the LP for the new top node of the parent offshoot. Since we have a warmstart basis from the top node of the new offshoot this can be done using a very good warmstart basis.

Obviously, the first two methods provide only a lower bound on the new optimal objective value. 

\subsection{Shortening Dives}

It is possible to implement the new method in a way that traditional branch-and-bound is just a special case. To this end we only need to ensure that in addition to storing open offshoots we can also store open nodes. This can be achieved, for example, by treating offshoots without a dive as normal nodes, which means when we select them we do not select an offshoot variable. Instead we treat it as the top node of a new offshoot directly. With this in place we can also have a limit on the number of bound changes in a dive. When the limit is hit, we store the last node as an open node in addition to storing the offshoot. The offshoot in this case does not end in a pruned node, but the method works regardless. If we set the limit of bound changes in a dive to zero, the method reverts to a traditional branch-and-bound method.

\subsection{Splitting offshoots}

Branching from the top on a bound change that was not the initial bound change of an offshoot invalidates all the internal objective values of the original nodes of an offshoot. This prevents us from pruning them and hurts the performance for very deep dives. Therefore it seems advantageous to split very long dives. However, this creates a non-terminal offshoot, so it extends our depth-first framework a little. For maximum efficiency, we need to resolve the linear relaxation of the bottom offshoot, so that we can have a valid objective value useful for pruning by bound. Which bound changes should be in the top of bottom half of the split is an interesting research question.

\section{Computational Evaluation}
The new method was prototyped using the MILP solver in SAS/OR. The prototype was meant as a way to evaluate the correctness and practicability of the method described and as such does not contain all the features and tricks of a full MILP solver. Nevertheless we present results using this prototype to give an impression of the capabilities of the new method.

The prototype plugs into the MILP solver after its root node when the actual branch-and-bound phase begins. It features a standard reliability branching strategy with a dynamic strong branching limit and a reliability limit of 5. For selecting the next offshoot we choose the best top bound first without explicit tie breaking. The prototype also features basic node presolve and reduced cost fixing techniques (only at the top of an offshoot), and also using the root reduced cost to fix columns globally as new incumbent solutions are found. What it notably lacks are more advanced node presolver techniques, local or global cuts in the nodes of the branch-and-bound tree, and primal heuristics.

We conduct our experiments on 96 machines running 2 jobs each on 16-core/2-socket Intel\textsuperscript{\tiny\textregistered} Xeon\textsuperscript{\tiny\textregistered} E5-2630 v3 @ 2.40GHz CPUs. All experiments are done with default settings, a memory limit of 62 GB, and a time limit of 2 hours. We use 798 instances that are the internal test set used to develop the SAS MILP solver. To evaluate our results we use performance profiles as described in \cite{Dolan2002}.

\subsection{Offshoot variable selection}

The first experiment evaluates several methods we implemented for choosing the offshoot variable, i.e., it is meant to judge the importance of shaping in the new method. We implemented four different methods:
\begin{description}
\item[bottom:] always branch from the bottom of an offshoot without changing the order of the variables;
\item[top:] always branch from the top of an offshoot without changing the order of the variables;
\item[pseudo:] choose the offshoot variable with the best reliable pseudocost score and branch from the top. If there are no offshoot variables with reliable pseudocost available, then choose the variable with the worst pseudocost score and branch from the bottom.
\item[pseudodual:] like \textbf{pseudo}, except if there are no variables with reliable pseudocost then use the variable with the worst dual information score and branch from the bottom. The dual information score is the reduced cost or the Farkas certificate of the pruned node when the offshoot was first processed.
\end{description}
The first two strategies do not shape the tree so they can be seen as a baseline for the performance of the method. The default method is \textbf{pseudodual}.

\begin{figure}
\includegraphics[width=\textwidth]{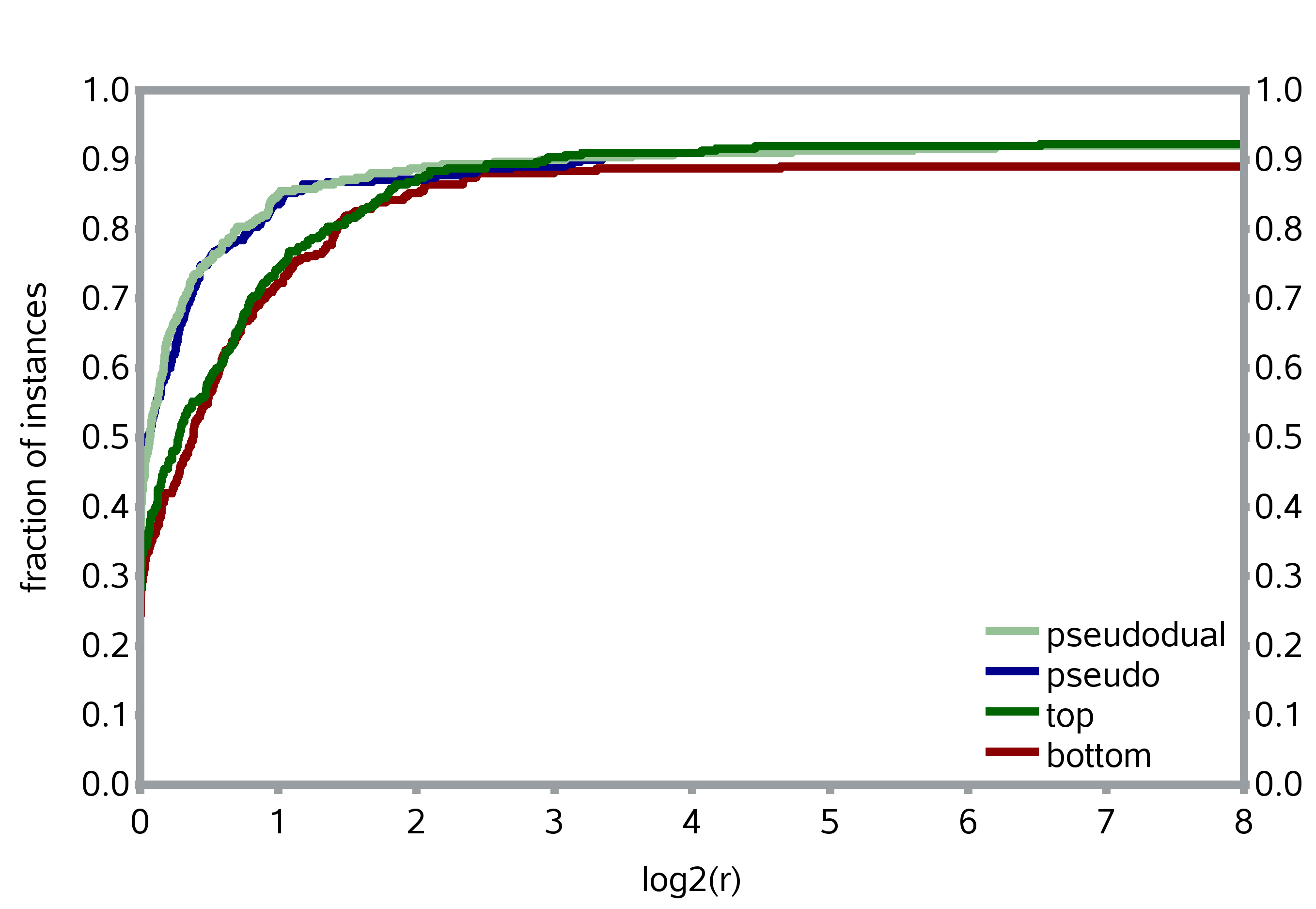}
\caption{Performance profile comparing offshoot variable selection strategies.}\label{shapingfig}
\end{figure}

Figure~\ref{shapingfig} shows the performance profile comparing the different offshoot variable selection strategies. In addition we would like to mention that the \textbf{pseudodual} strategy is about 13\% faster in the geometric mean of the solve times than the \textbf{bottom} strategy and solves 9 instances more within the time limit. We argue that this shows that shaping, at least in the context of this new method, has a clear impact on the performance. More advanced selection strategies can probably be developed that will demonstrate this even more profoundly.

\subsection{Trimming and pruning}

\begin{figure}
\includegraphics[width=\textwidth]{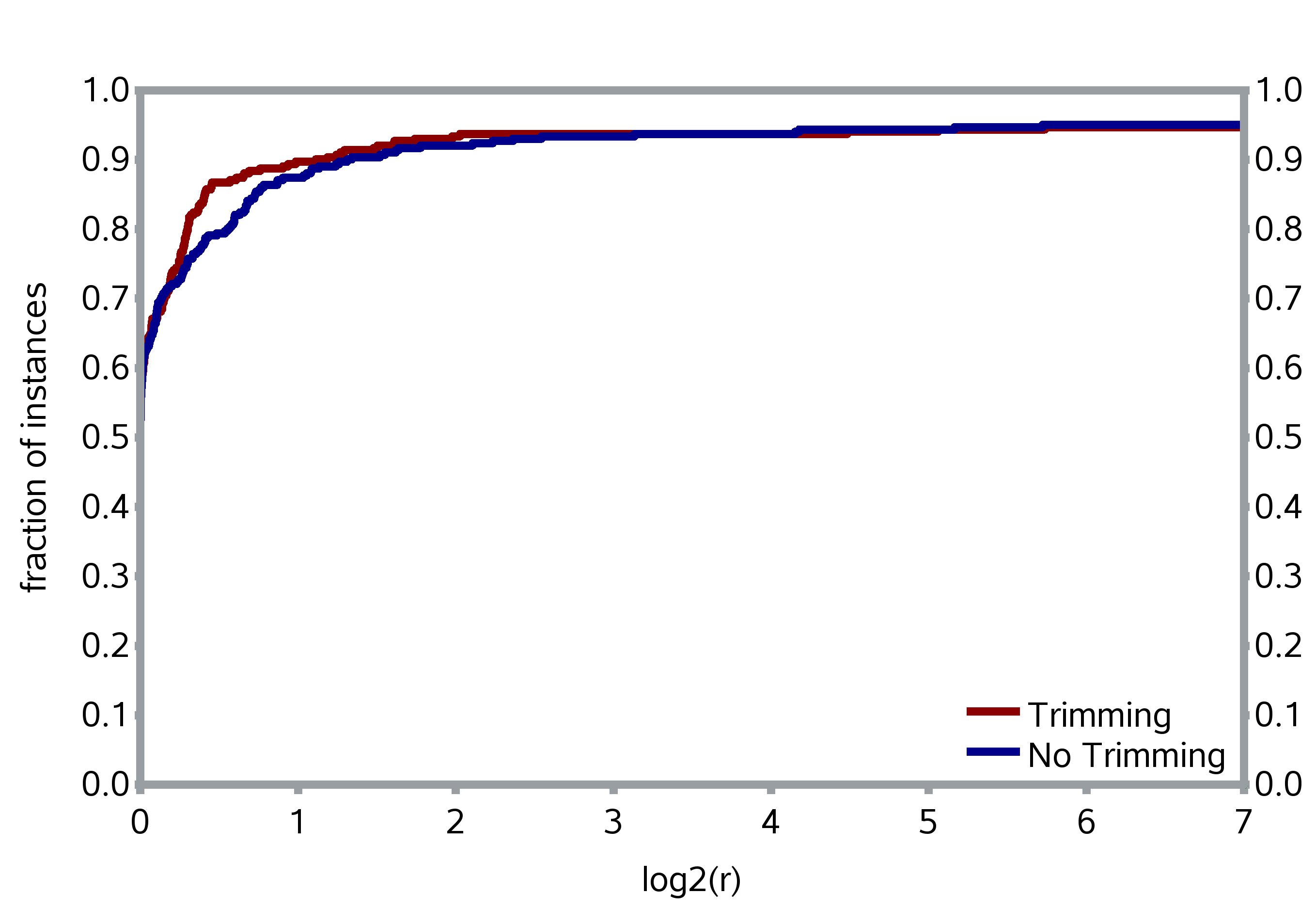}
\caption{Performance profile comparing a version of the prototype that does pruning and trimming with a version that does not.}\label{trimmingfig}
\end{figure}

Our second experiment is designed to show the combined importance of trimming and pruning. In our prototype implementation we delay trimming an offshoot until we need to choose an offshoot variable for the first time. Since we can either prune the bottom of the offshoot using the current cutoff or apply trimming using the dual information, we analyze how many reductions we get from either and choose the method that yields the most. In this experiment we compare the default version of our prototype that does this delayed pruning or trimming with a version where this feature has been disabled. The performance profile can be seen in Figure~\ref{trimmingfig}. The version with trimming and pruning is about 5\% faster in the geometric mean of solve times and solves 1 instance fewer within the time limit. Since this effect seems to be rather small we think that it would be necessary to look into better ways to trim dives. Some ideas are described in Section~\ref{advtrim}.

\subsection{Comparison against branch-and-bound}

\begin{figure}
\includegraphics[width=\textwidth]{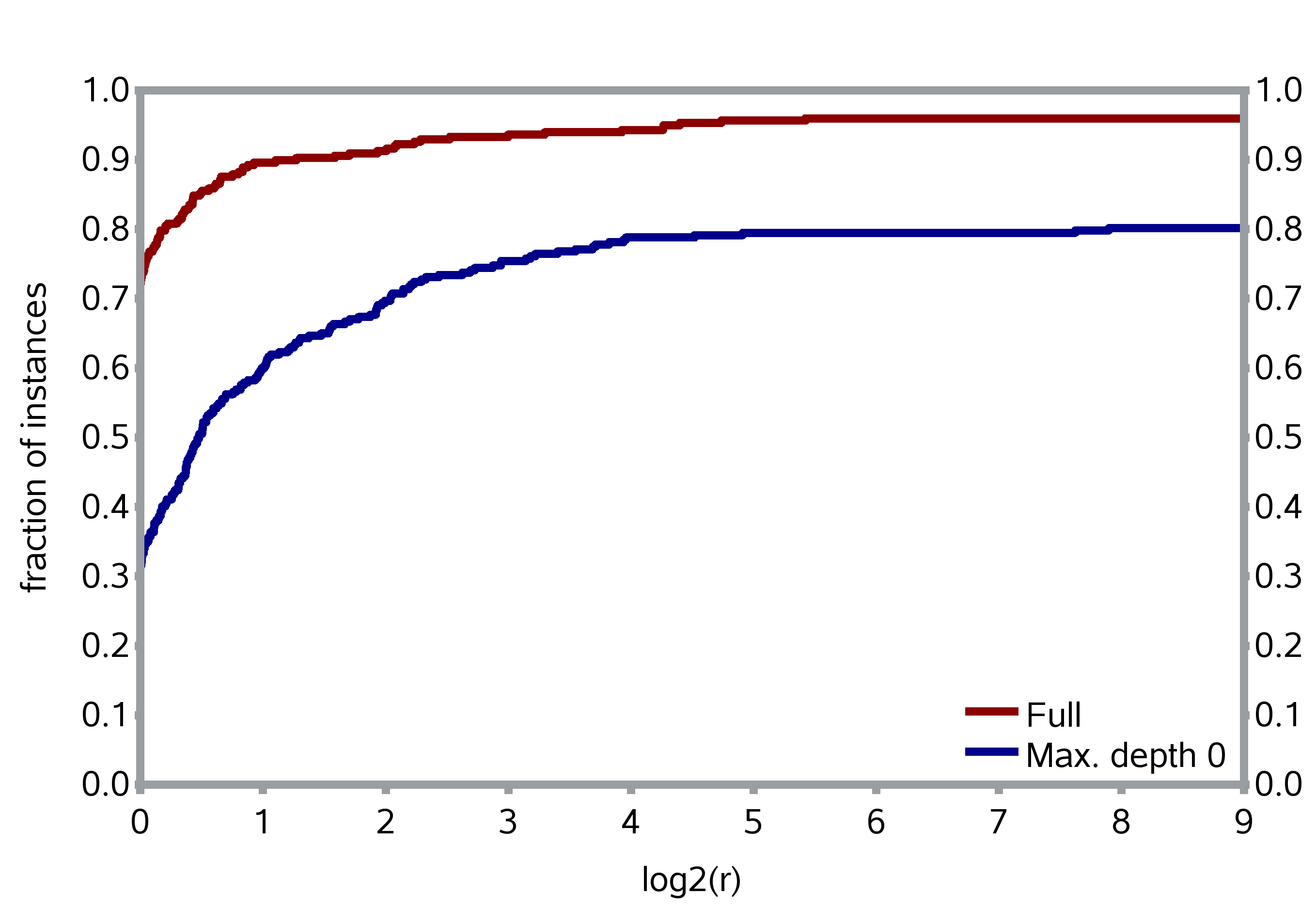}
\caption{Performance profile comparing the default version with version with a maximum depth of 0, i.e., that resembles a traditional branch-and-bound.}\label{depth0fig}
\end{figure}

In our final experiment we compare our default method that does not limit the depths of the dives to a version where the limit is 0. This means that the method with the limit is essentially a traditional branch-and-bound method. The comparison is not completely fair since a pure branch-and-bound method could be implemented more efficiently, especially regarding memory requirements. But it gives a first impression of how much could be gained by using our new method instead of a traditional branch-and-bound method. Figure~\ref{depth0fig} shows the performance profile. Our new method is 38\% faster in the geometric mean of the solve times and solves 47 instances more within the time limit. We consider this an encouraging result.

\section{Conclusions}
It will obviously take more research and a more elaborate implementation to see if our new method is superior to a traditional branch-and-bound method. From a theoretical perspective and from our preliminary experiments it seems likely that shaping and trimming the tree will result in improved run times. Even if the performance gains end up being very small there is also hope that our new method will result in a more stable performance.

So far we have not investigated other areas of application for our new method such as mixed integer non-linear optimization problems or branch-and-price algorithms. Since in these areas more flexibility in the tree might be even more advantageous we hope that it will find application there as well.

\providecommand{\bysame}{\leavevmode\hbox to3em{\hrulefill}\thinspace}
\providecommand{\MR}{\relax\ifhmode\unskip\space\fi MR }
\providecommand{\MRhref}[2]{%
  \href{http://www.ams.org/mathscinet-getitem?mr=#1}{#2}
}
\providecommand{\href}[2]{#2}


\begin{thebibliography}{1}

\bibitem{Achterberg2005}
T.~Achterberg, T.~Koch, and A.~Martin, \emph{Branching rules revisited},
  Operations Research Letters \textbf{33} (2005), no.~1, 42--54.

\bibitem{chinneck2008}
J.~W. Chinneck, \emph{Feasibility and infeasibility in optimization: Algorithms
  and computational methods}, International Series in Operations Research and
  Management Sciences, vol. 118, 2008.

\bibitem{chvatal1997}
V.~Chv{\'a}tal, \emph{Resolution search}, Discrete Applied Mathematics
  \textbf{73} (1997), no.~1, 81 -- 99.

\bibitem{Cook2012}
W.~Cook, \emph{Markowitz and {M}anne+ {E}astman+ {L}and and {D}oig= branch and
  bound}, Optimization Stories (2012), 227--238.

\bibitem{Dolan2002}
E.~Dolan and J.~Mor\'e, \emph{Benchmarking optimization software with
  performance profiles}, Mathematical Programming Series A \textbf{91} (2002),
  201--213.

\bibitem{Glover1976}
F.~Glover and L.~Tangedahl, \emph{Dynamic strategies for branch-and-bound},
  OMEGA - The International Journal of Management Science \textbf{4} (1976),
  no.~5, 571--576.

\bibitem{Hanafi2002}
S.~Hanafi and F.~Glover, \emph{Resolution search and dynamic branch-and-bound},
  Journal of Combinatorial Optimization \textbf{6} (2002), no.~4, 401--423.

\bibitem{Karzan2009}
F.~K{\i}l{\i}n{\c{c}}~Karzan, G.~Nemhauser, and M.~Savelsbergh,
  \emph{Information-based branching schemes for binary linear mixed integer
  problems}, Mathematical Programming Computation \textbf{1} (2009), no.~4,
  249--293.

\bibitem{Land1960}
A.~H. Land and A.~G. Doig, \emph{An automatic method of solving discrete
  programming problems}, Econometrica \textbf{28} (1960), no.~3, 497--520.

\end{thebibliography}
\end{document}